\documentclass[12pt,a4paper,reqno]{amsart}
\usepackage{amsmath}
\usepackage{amsfonts}
\usepackage{amssymb}
\numberwithin{equation}{section}

\theoremstyle{plain}

\newtheorem{theorem}[subsection]{Theorem}
\newtheorem{proposition}[subsection]{Proposition}
\newtheorem{lemma}[subsection]{Lemma}

\theoremstyle{definition}

\newtheorem{definition}[subsection]{Definition}

\newtheorem{question}[subsection]{Question}

\renewcommand{\leq}{\leqslant}
\renewcommand{\geq}{\geqslant}

\newsavebox{\proofbox}
\savebox{\proofbox}{\begin{picture}(7,7)%
  \put(0,0){\framebox(7,7){}}\end{picture}}
\def\boxeq{\tag*{\usebox{\proofbox}}}

\newcommand{\md}[1]{\ensuremath{(\mbox{mod}\, #1)}}

\newcommand\Spec{\operatorname{Spec}}
\newcommand\codim{\operatorname{codim}}
\newcommand\Supp{\operatorname{Supp}}
\def\proof{\noindent\textit{Proof. }}
\def\endproof{\hfill{\usebox{\proofbox}}}

\def\E{\mathbb{E}}
\def\Z{\mathbb{Z}}

\def\R{\mathbb{R}}
\def\F{\mathbb{F}}

\def\vs{\vspace{11pt}}

\begin{document}

\title{Boolean functions with small spectral norm}

\author{Ben Green}
\address{Department of Pure Mathematics and Mathematical Statistics\\
University of Cambridge\\
Wilberforce Road\\
Cambridge CB3 0WA\\
England } \email{b.j.green@dpmms.cam.ac.uk}

\author{Tom Sanders}
\address{Department of Pure Mathematics and Mathematical Statistics\\
University of Cambridge\\
Wilberforce Road\\
Cambridge CB3 0WA\\
England } \email{t.sanders@dpmms.cam.ac.uk}

\begin{abstract}
\noindent Let $f : \F_2^n \rightarrow \{0,1\}$ be a boolean
function, and suppose that the \emph{spectral norm} $\Vert f
\Vert_{A} := \sum_r |\widehat{f}(r)|$ of $f$ is at most $M$. Then
\[ f = \sum_{j = 1}^L \pm 1_{H_j},\] where $L \leq 2^{2^{CM^4}}$ and each $H_j$ is a subgroup of $\F_2^n$.\vspace{11pt}

\noindent This result may be regarded as a quantitative analogue
of the Cohen-Helson-Rudin structure theorem for idempotent
measures in locally compact abelian groups.

\end{abstract}

\thanks{The first author is a Clay Research Fellow, and thanks the Clay
Mathematics Institute for their support. Much of this work was
conducted while the second author was on a CMI-funded visit to
Boston, and he thanks the first author for arranging this and the
CMI for its support. Both authors would also like to thank the
Massachusetts Institute of Technology for their hospitality.}
\maketitle

\section{Introduction}

Let $G = \F_2^n$ be the $n$-dimensional cube, and let $f : G
\rightarrow \{0,1\}$ be a boolean function, or more generally a
function from $G$ to $\R$. In many works, particularly in
theoretical computer science, the Fourier transform
\[ \widehat{f}(r) := \E_{x \in G} f(x) (-1)^{r^T x} = \frac{1}{|G|}\sum_{x \in G} f(x) (-1)^{r^T x}\] is considered. Here, $r$ lies in the dual group $\widehat{G}$ which we have identified with $G$ by choosing the scalar product $u^Tv$.

It is natural to consider the $\ell^p$-norms
\[ \Vert \widehat{f} \Vert_p := \big( \sum_{r \in \widehat{G}} |\widehat{f}(r)|^p \big)^{1/p},\] for $1 \leq p < \infty$, as well as the $\ell^{\infty}$-norm
$\Vert \widehat{f} \Vert_{\infty} := \sup_r |\widehat{f}(r)|$.

There are many tools available for analysing these norms when $p
\geq 2$, particularly when $p$ is $\infty$ or an even integer.
When $1 \leq p < 2$, however, the situation is in many ways rather
mysterious. Of these cases, a very natural one is the endpoint $p
= 1$. In this case the norm $\Vert \widehat{f} \Vert_1$ is called
the \emph{algebra norm}, \emph{Wiener norm} or \emph{spectral
norm}; we shall denote it by $\Vert f \Vert_A$. It is quite easy
to show, using an instance of Young's inequality for convolutions,
that
\[ \Vert f g \Vert_A \leq \Vert f \Vert_A \Vert g \Vert_A\] for any two functions $f,g : G \rightarrow \R$. This explains the term \emph{algebra norm}.

A basic question is the following.

\begin{question}
Let $M$ be a fixed positive real number and let $f : G \rightarrow
\{0,1\}$ be a boolean function. When is $\| f\|_A \leq M$?
\end{question}

A partial answer to this question will be the main business of
this paper. By far the most important feature of the problem is
that we are asking it for \emph{boolean functions}, which take
only the values $0$ or $1$. There is a ready supply of functions
$f$ with $\|f\|_A$ small: Take for instance any pair $g,h:G
\rightarrow \mathbb{R}$. Then $f:=g \ast h$ has
\begin{equation*}
\|f\|_A = \|g \ast h \|_A = \|\widehat{g}\widehat{h}\|_1 \leq
\|\widehat{g}\|_2\|\widehat{h}\|_2 = \|g\|_2\|h\|_2,
\end{equation*}
which is small if $g$ and $h$ are small in $L^2$. It is rather
hard, however, to construct a large supply of such functions which
take only the values 0 and 1.

To get a feel for the question, we prove a simple folklore result
concerning the case $M = 1$. In fact, by choosing a suitable
argument from among the many available, one can cover the case $M
< 3/2$.

\begin{proposition}[Boolean functions with tiny spectral norm] Let $f : G \rightarrow \{0,1\}$ be a boolean function which does
not vanish identically. Then either $f = 1_{t + H}$, where $t + H$
is a coset of a subgroup of $G$, in which case $\Vert f \Vert_A =
1$, or else $\Vert f \Vert_A \geq 3/2$.
\end{proposition}
\proof First note that since $f$ is not identically zero we have
$\|f\|_\infty \geq 1$, and so $\|f\|_A \geq \|f\|_\infty \geq 1$,
by the simplest instance of the Hausdorff-Young inequality.

Now suppose that $H \leq G$. The Fourier transform of $1_{t + H}$
is supported on $H^{\perp} := \{r \in \widehat{G} : r^T x = 0 \; \mbox{for
all }x \in H\}$, and it has modulus $\Vert 1_H \Vert_1$ there. It
follows from this and the fact that $|H||H^{\perp}| = |G|$ that
$\Vert f \Vert_A = 1$ when $f = 1_{t + H}$.

To get the stronger statement claimed, we note that if $f$ is not
(the characteristic function of) a coset of a subgroup then there
are four distinct points $x,x+h,x+k,x+h+k$ forming a parallelogram
in $G$ such that $f(x)=f(x+h)=f(x+k)=1$ but $f(x + h + k)=0$ (this
is actually an if and only if statement -- we leave the proof of
both directions to the reader). Let \[ \phi := \delta_x +
\delta_{x + h} + \delta_{x + k} - \delta_{x + h + k},
\] thus $\phi(x) = \phi(x + h) = \phi(x + k) = |G|$, $\phi(x + h +
k) = -|G|$ and $\phi(y) = 0$ for all other $y$. Now we can compute
that
\[ \langle
f,\phi\rangle := \E_{x \in G} f(x) \phi(x) = 3 \textrm{ and }
\Vert \widehat{\phi} \Vert_{\infty} = 2,\] and so it follows from
Plancherel's theorem that
\[ 3 = \langle f,\phi\rangle = \langle \widehat{f},\widehat{\phi}\rangle \leq \Vert \widehat{f} \Vert_1 \Vert \widehat{\phi} \Vert_{\infty} \leq 2 \Vert f \Vert_A,\]
which proves the result.\endproof

\emph{Remarks.} We leave it to the reader to confirm that the
constant $3/2$ is best possible. The result (and proof) are
inspired by two papers of Saeki \cite{SSI,SSII} in which the
same question is addressed over all locally compact abelian groups
$G$. In that more general setting the constant $3/2$ should be
reduced to $\frac{1}{2}(1 + \sqrt{2})$, and equality can occur in
any group $G$ with an element of order $4$.

Returning to our main question, let us recall that $\Vert \cdot
\Vert_A$ is an algebra norm. Thus if $f_1,f_2 : G \rightarrow
\{0,1\}$ are functions for which $\Vert f_1 \Vert_A$ and $\Vert
f_2 \Vert_A$ are small then the functions $f_1 \vee f_2$, $f_1
\wedge f_2$,$1 - f_1$ and $1- f_2$ also have this property.
Loosely speaking, we refer to functions which can be obtained by a
small number of operations of this kind from the basic functions
$1_{t + H}$ as belonging to the \emph{coset ring} of $G$. In fact,
it is easy to see (ignoring quantitative issues for the time
being) that all elements of the coset ring are in fact of the form
\begin{equation}\label{eq0.3}
\sum_{j = 1}^L \pm 1_{H_j},\end{equation} for subgroups $H_j \leq
G$ and some ``small'' $L$.

One trivially has the bound
\[ \Vert \sum_{j = 1}^L \pm 1_{H_j} \Vert_A \leq L,\]
and so it is rather natural to ask whether something like the
converse is true; this is the main result of our paper.

\begin{theorem}[Main theorem]\label{mainthm}
Suppose that $f : G \rightarrow \{0,1\}$ has $\Vert f \Vert_A \leq M$. Then we may write
\[ f = \sum_{j = 1}^L \pm 1_{H_j},\] where the $H_j$ are subgroups of $G$ and $L \leq 2^{2^{CM^4}}$ for some absolute constant $C$.
\end{theorem}

\emph{Remarks.} The bound may seem unimpressive, and indeed in a
sense it is. However it depends only on $M$, a feature which we
believe is new to this paper. We do not dare to venture a guess as
to the correct bound, and it seems to us that it would be
difficult to use our method to reduce the number of exponentials
below two. It may be possible to reduce the power 4 somewhat,
although we have not attempted to do this.

The reader may wonder why we bothered to introduce the coset ring
at all, when only the very natural functions \eqref{eq0.3} are
involved in our theorem. The answer is that the description of the
coset ring in the form \eqref{eq0.3} is specific to the case $G =
\F_2^n$, and the phenomenon described by Theorem \ref{mainthm} is,
in a sense, more general.

Indeed our entire approach was motivated by Cohen's celebrated
idempotent theorem \cite{PJC}. Suppose that $G$ is a locally
compact abelian group, and that $M(G)$ is the Banach algebra of
finite measures on $G$ under convolution (see \cite[Appendix
E]{WR} for details). We say that a measure $\mu \in M(G)$
is \emph{idempotent} if $\mu \ast \mu = \mu$. Cohen's theorem is
that $\mu$ is idempotent if any only if $\widehat{\mu}$ lies in
the coset ring of $\Gamma = \widehat{G}$.

In our setting, Cohen's result implies that if $f : G \rightarrow \{0,1\}$ has $\|f\|_A < \infty$ then there is a
decomposition of the form of \eqref{eq0.3} with $L$ finite. This is, of course, a vacuous result. It was, however, natural to start with Cohen's argument (as described in
Rudin \cite{WR}) and try to make it effective. A
na\"{\i}ve attempt along these lines fails at several points and
there are even ``softer'' proofs of Cohen's theorem that we have
not managed to interpret in a finite setting at all, cf.
\cite{BHJFMFP}. Nevertheless access to these classical
results was crucial to our understanding and we could not have
written this paper without them.

We also import some ``modern'' ingredients from additive
combinatorics such as the Balog-Szemer\'edi theorem and Ruzsa's
analogue of Freiman's theorem. It seems to the authors that it may
be worth revisiting a number of classical results in the light of
these developments.

It is possible that our methods, in combination with the ideas in
\cite{PJC}, could lead to a fully quantitative proof of Cohen's
idempotent theorem. We intend to pursue this direction in future
work.

We conclude the introduction by remarking that the spectral norm
of boolean functions is discussed in the computer science
literature, but not in a great deal of detail. The papers
\cite{MB,YM} show that functions which can be computed
using a small binary decision tree have small spectral norm. Such
functions are, however, rather special elements of the coset ring.

\section{Notation}

Much of the notation we will use is implicit in the introduction,
but it may be helpful to clarify things here. When working with
functions on $G$, we will always use Haar probability measure.
Integration with respect to this measure will be denoted by $\E_{x
\in G}$, or sometimes just $\E$. If $f : G \rightarrow \R$ is a
function and $1 \leq p < \infty$ then we define
\[ \Vert f \Vert_p := \big( \E_{x \in G} |f(x)|^p \big)^{1/p}.\]
We also define $\Vert f \Vert_{\infty} := \sup_x |f(x)|$ as usual.
If $f_1,f_2 : G \rightarrow \R$ are two functions then we set
\[ \langle f_1 ,f_2 \rangle := \E_{x \in G} f_1(x)f_2(x)\]
and
\[ f_1 \ast f_2(x) := \E_{y \in G} f_1(x) f_2(x - y).\]
When working with the Fourier transforms of functions we will use
counting measure. Integration with respect to this measure will be
denoted by $\sum$ as usual. We defined the $\ell^p$ norms in the
introduction. Note also Plancherel's identity, which implies that
\[ \langle f_1,f_2\rangle = \langle \widehat{f}_1,\widehat{f}_2 \rangle := \sum_r \widehat{f}_1(r)\widehat{f}_2(r).\]
We will occasionally write, e.g., $E^{\wedge}$ when taking the
Fourier transform of a particularly complicated expression $E$.

Finally, a word concerning absolute constants. The letter $C$ will
always denote an absolute constant, but the exact value of this
constant may change form expression to expression. If in doubt,
the reader should recall that all instances of $C$ could, if
desired, be replaced by specific constants in such a way that all
our proofs are correct.

\section{Almost integer-valued functions and almost homomorphisms}

A key feature of this paper is that we cannot work entirely within
the ``category'' of boolean functions. We must also consider more
general functions which are close to being integer-valued.

\begin{definition}[Almost integer-valued functions]
Let $\epsilon \in (0,1/2)$. We say that a function $f : G
\rightarrow \R$ is \emph{$\epsilon$-almost integer-valued} if
there is a function $f_{\Z} : G \rightarrow \Z$ such that $\Vert f
- f_{\Z} \Vert_{\infty} \leq \epsilon$.
\end{definition}

We will need to study the behaviour of almost integer-valued
functions under a certain class of map. Let $H$ be a subgroup of
$G$. For any function $f : G \rightarrow \R$, we define $\psi_H f$
by
\[ (\psi_Hf)(x) := \E_{y \in x+H}f(y) = f \ast \mu_H(x),\]
where $\mu_H$ denotes the Haar probability measure on $H$.
Equivalently, one may define $\psi_H$ in terms of its Fourier
transform by
\[ (\psi_H f)^{\wedge}(r) := \widehat{f}(r)1_{H^{\perp}}(r),\]
where the subgroup $H^{\perp} \leq \widehat{G}$ is the annihilator
of $H$, defined by
\[ H^{\perp} := \{r \in \widehat{G} : r^T x = 0 \; \; \mbox{for all $x \in H$}\}.\]

The following simple properties of $\psi_H$ follow immediately
from the above definitions.

\begin{lemma}[Simple properties of $\psi_H$]\label{simple-prop}
The norm of $\psi_H$ is at most $1$ in both the operator norm
induced by the spectral norm and in that induced by the
$L^\infty$-norm. That is to say
\[ \Vert \psi_H f \Vert_A \leq \Vert f \Vert_A \;\textrm{ and }\; \Vert \psi_H f \Vert_{\infty} \leq \Vert f \Vert_{\infty} \;\textrm{ for all }\;f:G \rightarrow \mathbb{R}.\]
\end{lemma}

\begin{definition}[Spectral support]\label{def3.1}
Let $\eta > 0$ be a parameter, let $f : G \rightarrow \R$ be a
function, and suppose that $H \leq G$. Then we say that $f$ is
\emph{$\eta$-spectrally supported} on $H$ if
\[ \sup_{ r \not\in H^\perp}{\sum_{r' \in r + H^{\perp}} |\widehat{f}(r')|} \leq \eta.\]
\end{definition}

Note that we do \emph{not} assume that $\widehat{f}$ has
substantial mass on $H^{\perp}$ itself.

\begin{lemma}[Finding the spectral support]\label{lem3.2}
Let $H \leq G$ be any subgroup, let $\eta > 0$ be any parameter,
and let $f : G \rightarrow \R$ be a function with $\Vert f \Vert_A
\leq M$. Then there is a subgroup $H' \leq H$ with
\[ \codim(H: H') \leq M/\eta\]
such that $f$ is $\eta$-spectrally supported on $H'$.
\end{lemma}
\proof Set $H_0 := H$. We define a descending sequence $H_0 \geq
H_1 \geq \dots$ of subgroups with $\codim(H_i : H_{i+1}) = 1$.

If, at some stage, $f$ is $\eta$-spectrally supported on $H_i$
then we stop. If not, there is some $r_i \notin H_i^{\perp}$ such
that
\[ \sum_{r' \in r_i + H_i^{\perp}} |\widehat{f}(r')| > \eta.\]
Define $H_{i+1}^{\perp}$ to be the subgroup of $G$ generated by
$r_i$ and $H_i^{\perp}$. It is clear that for any $j$ we have
\[ \Vert f \Vert_A \geq \sum_{i = 0}^{j} \sum_{r' \in r_{i} + H_i^{\perp}} |\widehat{f}(r')|,\] and so this inductive process must terminate after no more than $M/\eta$ steps.\endproof

The purpose of Definition \ref{def3.1} and Lemma \ref{lem3.2} is to allow us to use the following approximate homomorphism property.

\begin{lemma}[$\psi_H$ is an approximate homomorphism]\label{lem3.3}
Suppose that $f,g : G \rightarrow \R$ are two functions, and that $f$ is $\eta$-spectrally supported on $H$. Then
\[ \Vert \psi_H(fg) - \psi_H(f)\psi_H(g) \Vert_A \leq \eta \Vert g \Vert_A.\]
\end{lemma}
\proof We have
\[ \Vert \psi_H(fg) - \psi_H(f)\psi_H(g)\Vert_A = \sum_r \big| \widehat{f}\ast \widehat{g}(r)1_{H^{\perp}}(r) - \widehat{f}1_{H^{\perp}} \ast \widehat{g}1_{H^{\perp}}(r) \big|.\]
However
\begin{align*}
\widehat{f}\ast \widehat{g}(r)1_{H^{\perp}}(r) &= \sum_s
\widehat{f}(r-s) \widehat{g}(s)1_{H^{\perp}}(r) \\ &= \sum_s
\widehat{f}(r-s) \widehat{g}(s)1_{H^{\perp}}(r-s)1_{H^{\perp}}(s)
+ \sum_{s \notin H^{\perp}} \widehat{f}(r-s)
\widehat{g}(s)1_{H^{\perp}}(r) \\ &= \widehat{f}1_{H^{\perp}} \ast
\widehat{g}1_{H^{\perp}}(r) + \sum_{s \notin H^{\perp}}
\widehat{f}(r-s) \widehat{g}(s)1_{H^{\perp}}(r).
\end{align*}
Now simply note that
\begin{align*} \sum_r \big| \sum_{s \notin H^{\perp}} \widehat{f}(r-s) \widehat{g}(s)1_{H^{\perp}}(r) \big| &\leq \sum_{s \notin H^{\perp}} \sum_{r \in H^{\perp}} |\widehat{f}(r-s)||\widehat{g}(s)| \\ &\leq \sum_{s}|\widehat{g}(s)| \sup_{s \notin H^{\perp}} \sum_{r\in H^{\perp}} |\widehat{f}(r - s)| \\ &\leq \eta \Vert g \Vert_A.\end{align*}
This completes the proof.\endproof

Our aim is to show that, provided the parameter $\eta$ is suitably small, the map $\psi_H$ preserves almost integer-valued functions.

\begin{lemma}\label{poly-lem}
Let $d \geq 0$ be an integer, and write $P_d(X) :=
4^d(2d)!^{-1}\prod_{j = -d}^d (X - j)$. Let $\epsilon,\delta \leq
1/2$ be positive real parameters. Let $f : G \rightarrow \R$ be a
function.
\begin{enumerate}
\item If $f$ is $\epsilon$-almost integer-valued and $\Vert f
\Vert_{\infty} \leq d$ then $\Vert P_d(f) \Vert_{\infty} \leq
\epsilon 4^d$. \item If $\Vert P_d(f) \Vert_{\infty} \leq \delta$
then $f$ is $\delta$-almost integer-valued.
\end{enumerate}
\end{lemma}
\proof To prove the first statement, simply note that
\[ \Vert P_d (f) \Vert_{\infty} \leq 4^d(2d)!^{-1} \cdot (2d)! \sup_x \inf_j |f(x) - j| \leq  \epsilon 4^d. \]
To prove the second, observe that
\[ |P_d(t)| \geq |P_d(\overline{t})|\] for all $t \in \R$, where $\overline{t} \equiv t \md{1}$ and lies in the interval $(-1/2,1/2]$. Furthermore one may easily confirm that
\[ |P_d(\overline{t})|  \geq \frac{|\overline{t}| \cdot 4^dd! \prod_{j = 0}^{d-1} (j + \textstyle\frac{1}{2})}{(2d)!} \geq |\overline{t}| \] for all $t$.
It follows that for all $x$ the distance from $f(x)$ to the nearest integer is no more than $\delta$.\endproof

\emph{Remark.} We have normalised the polynomials $P_d$ slightly
arbitrarily, so that no factors were lost in (2). This makes no
essential difference to the argument.

To apply this, we combine it with the rest of the results of this
section to obtain the following corollary.

\begin{proposition}[$\psi_H$ preserves almost integer-valued functions]\label{prop3.4} Suppose that $f : G \rightarrow \R$ is $\epsilon$-almost
integer-valued, and that $\Vert f \Vert_A \leq M$ for some $M \geq
1/2$. Suppose that $\eta \leq 2^{-CM(1+\log M)}\epsilon$ for
some suitably large $C$ and that $f$ is $\eta$-spectrally
supported on $H$. Then both $\psi_H f$ and $f - \psi_H f$ are
$(2^{CM}\epsilon)$-almost integer-valued.\end{proposition} \proof
It clearly suffices to prove the result for $\psi_H f$, as the sum
or difference of two almost integer-valued functions is almost
integer-valued. Set $d := \lceil M \rceil$. Since $\Vert f
\Vert_{\infty} \leq \Vert f \Vert_A \leq M$, Lemma \ref{poly-lem}
(1) implies that
\[ \Vert P_d (f) \Vert_{\infty} \leq 2^{CM}\epsilon.\]
From Lemma \ref{simple-prop} it follows that
\begin{equation}\label{eq3-27} \Vert \psi_H(P_d(f)) \Vert_{\infty} \leq 2^{CM}\epsilon.  \end{equation}
Now an easy induction based on Lemma \ref{lem3.3} (and Lemma
\ref{simple-prop}) confirms that
\[ \Vert \psi_H(f^n) - (\psi_H f)^n \Vert_A \leq \eta (n-1)M^{n-1},\] for any positive integer $n$. It follows from this that
\[ \Vert \psi_H(P_d(f)) - P_d(\psi_H f) \Vert_A \leq 2^{CM(1 + \log M)}\eta \leq 2^{CM}\epsilon ,\]
and hence in view of \eqref{eq3-27} that
\[ \Vert P_d(\psi_H f) \Vert_{\infty} \leq 2^{CM}\epsilon.\]
The result is now an immediate consequence of Lemma \ref{poly-lem} (2).\endproof

\section{Ruzsa's analogue of Freiman's theorem}\label{sec4}

In the next two sections we use variants of a well-known sequence
of arguments in additive combinatorics. The objective is to prove Proposition \ref{concentrate-prop}, which roughly speaking states that a function with small $A(G)$-norm concentrates on a subspace. We will supply
original references for the results we use, but would also
recommend the book \cite{TCTVHV} as a general resource for
this subject.

When we actually prove Proposition \ref{concentrate-prop} we will
find ourselves dealing with a set $A \subseteq G$ with \emph{small
doubling}, that is to say a set $A$ with $\E 1_{A + A} \leq K\E
1_A$ for some ``not too large'' $K$. There is a beautiful theorem
of Imre Ruzsa \cite{IZRArb} (see also \cite{TSFT})
which states that in this case $A$ is contained in a subgroup $H
\leq G$ with density at most $K^2 2^{K^4}\E 1_A$. This certainly implies that
\begin{equation}\label{eq3-48} \E_{x \in A}{1_H(x)}\geq 1 \;\textrm{ and }\; \E_{x \in H}{1_A(x)}\geq 2^{-CK^C}.\end{equation}
One could use this result as it is, and obtain a bound in Theorem
\ref{mainthm} with a three-fold iterated exponential. To reduce
the number of exponentials to two, we need a different version of
Ruzsa's result, in which we shall replace \eqref{eq3-48} with
\begin{equation*} \E_{x \in A}{1_H(x)}\geq 2^{-CK^C} \textrm{ and } \E_{x \in H}{1_A(x)}\geq cK^{-C}.\end{equation*}
A more precise version of the following proposition, which is the main result
of the section, will be contained in a forthcoming paper of the
first author and Terence Tao. The authors are grateful to the latter
for useful discussions regarding this circle of ideas.

\begin{proposition}[Freiman in torsion groups, refined]\label{ruz-refinement}
Suppose that $A \subseteq G$ is a set with $\E 1_{A+A} \leq K\E
1_A$. Then there is a subgroup $H \leq G$ such that
\[ \E_{x \in A}{1_H(x)}\geq 2^{-CK^C} \;\textrm{ and }\; \E_{x \in H}{1_A(x)}\geq cK^{-C}.\]
\end{proposition}

\emph{Remark.} It is an important unsolved problem to decide
whether or not one may replace $2^{-CK^C}$ by a polynomial in
the first bound. This is known as the Polynomial Freiman-Ruzsa
conjecture (PFR); see for example \cite{BJGFFM}. The
truth of this conjecture, however, would not make an essential
difference to the bound we obtain in Theorem \ref{mainthm}.

To prove Proposition \ref{ruz-refinement} we need to set up a
little notation. Write $\alpha := \E 1_A$ for the density of $A$ in $G$. Put
\[ \nu^{(4)}(x) := 1_A \ast 1_A \ast 1_{A} \ast 1_{A} (x),\] and
for any parameter $\eta > 0$ define
\[ S_{\eta} := \{x \in G : \nu^{(4)}(x) \geq \eta \alpha^3\}.\]
For a parameter $\rho \in (0,1)$, we write (as is becoming
standard)
\[ \Spec_{\rho}(A) := \{r \in \widehat{G} : |\widehat{1}_A(r)| \geq \rho \alpha\}.\]
We begin by recording a well-known argument of Bogolyubov \cite{NNB} in this
language.
\begin{lemma}[Bogolyubov's argument]\label{bog} Suppose that $A \subset G$. Let $\delta,\epsilon \in (0,1)$ be any parameters and set $\rho :=
(\epsilon/2)^{1/2}$ and $H:=\Spec_{\rho}(A)^\perp$. Then
\[ S_{\delta} + H \subseteq S_{\delta - \epsilon}.\]
\end{lemma}
\proof Suppose that $x \in S_{\delta}$ and that $h \in H$. Then we have
\begin{align*}
\nu^{(4)}(x + h) &= \sum_{r \in \widehat{G}} |\widehat{1}_A(r)|^4 (-1)^{r^T (x + h)} \\ &= \sum_{r \in \Spec_{\rho}(A)} |\widehat{1}_A(r)|^4 (-1)^{r^T x} + \sum_{r \notin \Spec_{\rho}(A)} |\widehat{1}_A(r)|^4 (-1)^{r^T (x + h)} \\ & \geq \nu^{(4)}(x) - 2\sum_{r \notin \Spec_{\rho}(A)} |\widehat{1}_A(r)|^4.
\end{align*}
Thus we only need observe, using Parseval's identity and the definition of $\Spec_{\rho}(A)$, that
\begin{equation}\boxeq
\sum_{r \notin \Spec_{\rho}(A)} |\widehat{1}_A(r)|^4 \leq \sup_{r
\not \in \Spec_\rho(A)}{|\widehat{1}_A(r)|^2} \sum_{r \in
G}{|\widehat{1}_A(r)|^2} \leq (\rho\alpha)^2\cdot \alpha =
\epsilon \alpha^3/2.
\end{equation}

The next two lemmas are the vehicles by which we leverage the
assumption that $A$ has small doubling. The first states that $A$ has large density on a translate of $S_{\eta}$, provided $\eta$ is sufficiently small.

\begin{lemma}\label{lem1.3}
Suppose that $A \subseteq G$, that $\E 1_{A+A} \leq K \E 1_A$ and that $\eta \leq 1/2K^4$ is a parameter. Then $\E 1_{S_{\eta}} \geq \alpha/2$ and
\begin{equation*} \Vert 1_A \ast 1_{S_\eta}\Vert_{\infty} \geq
\eta\alpha/2.
\end{equation*}
\end{lemma}
\proof Averaging $\nu^{(4)}$ over $x \in G$, we get
\begin{align*} \alpha^4 & = \E \nu^{(4)} \\ & \leq \alpha^3\E 1_{S_{\eta}}  + \eta\alpha^3 \E 1_{4A} \\ &\leq \alpha^3\E 1_{S_{\eta}}  + \eta\alpha^4 K^4 \\ &\leq \alpha^3\E 1_{S_{\eta}}  + \alpha^4/2,\end{align*}
where the second inequality follows from the Pl\"unnecke--Ruzsa
inequalities \cite{IZR}, and the third from the condition on
$\eta$. The first conclusion of the lemma follows immediately upon rearranging.

For the second part we use the first to see that
\begin{align*}
\eta \alpha^4/2  &\leq \eta \alpha^3 \E 1_{S_\eta}  \\ &\leq \E (1_{S_\eta} \nu^{(4)}) \\ & =  \langle 1_{S_\eta},1_A
\ast 1_A \ast 1_A \ast 1_A \rangle\\ & =  \langle 1_A \ast
1_{S_\eta}, 1_A \ast 1_A \ast 1_A\rangle\\ & \leq  \Vert 1_A \ast 1_{S_{\eta}} \Vert_{\infty} \E 1_A \ast 1_A \ast 1_A.
\end{align*}
The conclusion follows immediately since $\E 1_A \ast 1_A \ast 1_A =
\alpha^3$.\endproof\vs

\begin{lemma}\label{lem1.4} Suppose that $A \subseteq G$ has density $\alpha := \E 1_A$ and that $\E
1_{A+A} \leq K\E 1_A$. Then there is a subgroup $H \leq G$ with
$\E 1_H \geq (\alpha/2)^{-CK^{12}}$ and some $\eta$ with $1/4K^4
\leq \eta \leq 1/2K^4$, such that $S_{\eta}$ is a union of cosets
of $H$ together with an exceptional set $X$ satisfying $\E 1_X \leq
\alpha/16K^4$.
\end{lemma}
\proof Let $\eta_0 := 1/2K^4$ and set $\epsilon := 1/64K^{12}$. Consider the nested sequence
\[ S_{\eta_0} \subseteq S_{\eta_0 - \epsilon} \subseteq \dots \subseteq S_{\eta_0 - (L-1)\epsilon},\] where $L := 1/4K^4\epsilon$. By the Pl\"unnecke-Ruzsa inequalities we have
\cite{IZR}
\[ \E 1_{S_{\eta}} \leq \E 1_{4A} \leq K^4 \alpha\] for any $\eta$, and therefore by the pigeonhole
principle there is some $j$ with $0 \leq j < L$ such that
\[ \E 1_{S_{\eta_0 - (j+1)\epsilon} \setminus S_{\eta_0 - j\epsilon}} \leq \alpha/16K^4.\]

Now we apply Lemma \ref{bog}. Writing $H = \Spec_{\rho}(A)^{\perp}$ where $\rho := 1/16K^6$, we know from that lemma that
\[ S_{\eta_0 - j\epsilon} + H \subseteq S_{\eta_0 - (j+1)\epsilon}.\]
Thus $S_{\eta_0 - (j+1)\epsilon}$ can be written as a union of
cosets of $H$ together with an exceptional set $X$ of density at
most $\alpha/16K^4$.

It remains to establish the claimed lower bound on $\E 1_H$. By a
lemma of Chang \cite[Lemma 3.1]{MCC} (see also
\cite[Lecture 14, Lemma 3]{BJGRKP} and \cite{WR.TSG,TCTVHV}) the set
$\Spec_{\rho}(A)$ is contained in a subgroup of $G$ with dimension
$O( \rho^{-2}(1+\log(1/\alpha))) = O(K^{12}\log \alpha^{-1})$.
This concludes the proof.\endproof\vs

\emph{Proof of Proposition \ref{ruz-refinement}.} It is sufficient
to prove the proposition when $G=\langle A \rangle$ in which case, by
Ruzsa's Theorem \cite{IZRArb}, we have $\alpha \geq
K^{-2}2^{-K^4}$.

Apply Lemma \ref{lem1.4} to get a subgroup $H'$ with $\E 1_{H'}
\geq 2^{-CK^{16}}$ and some $\eta$ with $1/4K^4 \leq \eta \leq
1/2K^4$, such that
\[ S_\eta = X \cup \bigcup_{y \in Y} (y + H'),\]
where $\E 1_X \leq \alpha/16K^4$. Writing $\mu_{H'} := 1_{H'}/\E1_{H'}$ for the Haar measure on $H'$, we have for all $x \in G$ that
\begin{align*}
1_A \ast 1_{S_{\eta}} \ast \mu_{H'}(x) &\geq 1_A \ast 1_{S_{\eta}\setminus X} \ast \mu_{H'} (x) \\ & = 1_A \ast 1_{S_{\eta} \setminus X}(x) \\ & = 1_A \ast 1_{S_{\eta}}(x) - 1_A \ast 1_X (x) \\ & \geq 1_A \ast 1_{S_{\eta}}(x) - \E 1_X \\ & \geq 1_A \ast 1_{S_{\eta}}(x) - \alpha/16K^4.
\end{align*}
It follows from this, Lemma \ref{lem1.3} and the assumption that $\eta \geq 1/4K^4$ that
\begin{equation}\label{use-soon} \Vert 1_A \ast 1_{S_{\eta}} \ast \mu_{H'} \Vert_{\infty} \geq \Vert 1_A \ast 1_{S_{\eta}} \Vert_{\infty} - \alpha/16K^4 \geq \eta \alpha/2 - \alpha/16K^4 \geq \alpha/16K^4.
\end{equation}
Furthermore by the Pl\"unnecke-Ruzsa inequalities \cite{IZR} we have
\[
\Vert 1_A \ast 1_{S_{\eta}} \ast \mu_{H'}\Vert_{\infty}  \leq \Vert 1_A \ast \mu_{H'}\Vert_{\infty} \E 1_{S_{\eta}}  \leq \Vert 1_A \ast \mu_{H'} \Vert_{\infty} \E1_{4A}  \leq K^4 \alpha \Vert 1_A \ast \mu_{H'} \Vert_{\infty}.\]
Comparing with \eqref{use-soon} leads immediately to
\[ \Vert 1_A \ast \mu_{H'} \Vert_{\infty} \geq 1/16K^8.\]
We have found a coset of $H'$ on which the relative density of $A$
is at least $1/16K^8$; by adjoining the zero element to $H'$ if
necessary, one obtains a subgroup $H$ on which the relative
density of $A$ is at least $1/32K^8$, that is to say
\[ \E_{x \in H} 1_A(x) \geq 1/32K^8.\]
To complete the proof of Proposition \ref{ruz-refinement} it remains to note that
\begin{equation}\boxeq
\E_{x \in A}1_{H}(x) = \frac{\E 1_H}{\E 1_A}\cdot \E_{x \in
H}{1_A(x)} \geq (\E 1_H) \E_{x \in H}{1_A(x)} \geq 2^{-CK^{16}}.
\end{equation}

\section{Concentration on a subgroup}

\begin{proposition}[Concentration on a subgroup]\label{concentrate-prop}
Suppose that $f : G \rightarrow \R$ is an
$\epsilon$-integer-valued function with $\Vert f \Vert_A \leq M$,
where $M \geq 1/2$ and $\epsilon \leq 2^{-C M^4}$. Then there is
a subgroup $H \leq G$ with
\begin{equation*}
\E 1_H \geq 2^{-2^{CM^4}}\Vert f_{\Z}\Vert_1
\end{equation*}
and
\[  \sup_{x
\in G}{|\E_{y \in x+H}{f(y)}|} = \|\psi_H f \|_\infty \geq
2^{-CM^4}.\]
\end{proposition}

\begin{definition}[Arithmetic connectedness]\label{arith-connect-def}
Let $m$ be a positive integer. Suppose that $A \subseteq G$ is a
set with $0 \notin A$. Then we say that $A$ is
\emph{$m$-arithmetically connected} if, for any choice of distinct
$a_1,\dots,a_{m} \in A$, one of the following alternatives holds:
\begin{enumerate}
\item The vectors $a_1, \dots, a_m$ are linearly dependent;
\item The vectors $a_1, \dots, a_m$ are linearly independent but there exists a further $a' \in A$ such that $a'$ lies in the linear span of the $a_i$.
\end{enumerate}
\end{definition}

The next lemma imports the tools we developed in \S\ref{sec4}. The
result allows us to weaken the condition of small doubling in
Proposition \ref{ruz-refinement} to that of arithmetic
connectedness.

\begin{lemma}\label{additive-comb-lem}
Suppose that $m$ is a positive integer and that $A \subseteq G$ is
a set with $0 \notin A$. Suppose that $A$ is $m$-arithmetically
connected. Then there is a subgroup $H \leq G$ such that
\begin{equation*}
\E_{x \in A}{1_H(x)} \geq 2^{-2^{Cm}} \;\textrm{ and } \; \E_{x \in
H}{1_A(x)} \geq 2^{-Cm}.
\end{equation*}
\end{lemma}
\proof If $|A| < m^2$ the result is trivial, so we stipulate that $|A| \geq m^2$. Pick any $m$-tuple $(a_1,\dots,a_m)$ of distinct elements of $A$. With the stipulated lower bound on $|A|$, there are at least $|A|^m/2$ such $m$-tuples. We know that either the vectors $a_1,\dots,a_m$ are linearly dependent, or else there is a further $a' \in A$ such that $a'$ lies in the linear span of the $a_i$. In either situation there is some linear relation
\[ \lambda_1 a_1 + \dots + \lambda_m a_m + \lambda' a' = 0\]
where $\vec{\lambda} := (\lambda_1,\dots,\lambda_m,\lambda')$ has
elements in $\F_2$ and, since $0 \notin A$ and the $a_i$s are
distinct, at least three of the components of $\vec{\lambda}$ are
nonzero. By the pigeonhole principle, it follows that there is
some $\vec{\lambda}$ such that the linear equation
\[ \lambda_1 x_1 + \dots + \lambda_m x_m + \lambda' x' = 0\] has at least $|A|^m/2^{m+2}$ solutions with $x_1,\dots,x_m,x' \in A$. Removing the zero coefficients, we may thus assert that there is some $r$, $3 \leq r \leq m + 1$, such that the equation
\[ x_1 + \dots + x_r = 0\]
has at least $|A|^{r-1}/2^{m+2}$ solutions with $x_1,\dots,x_r \in A$. Note that this is a strong structural statement about $A$, since the maximum possible number of solutions to such an equation is $|A|^{r-1}$.

We claim that there are at least $2^{-2m-4}|A|^3$ solutions to
$x_1 + x_2 = x_3 + x_4$ with $x_i \in A$. To see this, write
$R_l(x)$ for the number of $l$-tuples $(x_1,\dots,x_l) \in A^l$
such that $x_1 + \dots + x_l = x$, and note that
\[ \sum_x R_2(x)R_{r-2}(x) =R_r(0) \geq |A|^{r-1}/2^{m+2}.\]
Noting that $R_{r-2}(x) \leq |A|^{r-3}$ for all $x$ (here, of course, it is important that $r \geq 3$) we see from the Cauchy-Schwarz inequality that
\[ \sum_x R_2(x)^2 \geq \frac{|A|^{2(r-1)}}{2^{2m+4} \big( \sum_x R_{r-2}(x)^2 \big)} \geq 2^{-2m-4}|A|^3,\] confirming the claim.

It now follows from the Balog-Szemer\'edi-Gowers theorem
\cite[Proposition 12]{WTG} that there is some set $A'
\subseteq A$, $\E 1_{A'} \geq 2^{-Cm}\E 1_A$, such that $\E1_{A' + A'} \leq
2^{Cm}\E1_{A'}$. By Proposition \ref{ruz-refinement} there is a
subgroup $H \leq G$ with
\[ \E_{x \in A'}{1_H(x)} \geq  2^{-2^{Cm}} \; \textrm{ and } \; \E_{x \in H}{1_{A'}(x)} \geq 2^{-Cm}.\]
The result follows since $\E 1_A \geq \E 1_{A'} \geq 2^{-Cm}\E
1_A$.\endproof

\emph{Proof of Proposition \ref{concentrate-prop}.} We begin by
decomposing $f^2$ as $g+h$ where $g=f_{\Z}^2$ and $h=f^2-g$. We
have
\[ \Vert h \Vert_{\infty} = \Vert f^2 - f_{\Z}^2 \Vert_{\infty} \leq \Vert f - f_{\Z} \Vert_{\infty} \Vert f + f_{\Z} \Vert_{\infty} \leq 4\epsilon M,\]
the latter inequality being a consequence of the fact that $\Vert f \Vert_{\infty} \leq \Vert f \Vert_A \leq M$ and that $\Vert f_{\Z} \Vert_{\infty} \leq \Vert f \Vert_{\infty} + \epsilon$.

Set $m = \lceil (2M)^4 \rceil$ and suppose that $V \leq G$ is a
subgroup of dimension $m$. We have
\[ \Vert f^2 1_{V} \Vert_A \leq \Vert f \Vert_A^2 \Vert 1_{V} \Vert_A = \Vert f \Vert_A^2  \leq M^2.\]
In view of the trivial estimate
\[ \Vert h 1_{V} \Vert_{A} \leq \sum_{x \in V}  \Vert h(x) 1_x \Vert_A \leq 2^m \Vert h \Vert_{\infty} \leq 2^{m+2}\epsilon M,\] it follows from the triangle inequality and the assumption on
$\epsilon$ that
\begin{equation}\label{g-l1-bd} \Vert g 1_{V} \Vert_A \leq M^2 + 2^{m+2}\epsilon M \leq 2M^2.\end{equation}

Write $A := \Supp(g) = \Supp(f_{\Z})$. If $A$ is all of $G$ then
the proposition follows trivially so we may assume that this is
not the case. Hence by replacing $f(x)$ by $f(x + y)$ for some $y
\notin A$, we may assume without loss of generality that $0 \notin
A$. We claim that $A$ is $m$-arithmetically connected in the sense
of Definition \ref{arith-connect-def}. If this is not the case
then there are elements $a_1,\dots,a_m \in A$ such that the
vectors $a_1,\dots,a_m$ are linearly independent, and such that
there is no $a' \in A$ with $a'$ in the linear span of the $a_i$.
Writing $V$ for the subgroup of $G$ spanned by the $a_i$, this
means that the support of $g 1_{V}$ is precisely
$\{a_1,\dots,a_m\}$.

Thus we have
\[ g1_V(x) = \sum_{i = 1}^m g(a_i) 1_{a_i}(x).\]
Now we may compute that
\[ \Vert (g1_V)^{\wedge} \Vert_2^2 = \Vert g1_V \Vert_2^2 = \frac{1}{|G|}\sum_{i=1}^m |g(a_i)|^2 \geq \frac{m}{|G|}\]
and that
\[ \Vert (g1_V)^{\wedge} \Vert_4^4 = \frac{1}{|G|^3} \sum_{\substack{i_1,i_2,i_3,i_4 \\ a_{i_1} + a_{i_2} = a_{i_3} + a_{i_4}}} |g(a_i)|^4 \leq \frac{3}{|G|^3} \big( \sum_{i = 1}^m |g(a_i)|^2 \big)^2 = \frac{3}{|G|} \Vert (g1_V)^{\wedge} \Vert_2^4,\] the middle inequality following from the observation that $a_{i_1}
+ a_{i_2} = a_{i_3} + a_{i_4}$ only if $i_1=i_2$ and $i_3=i_4$, or
$i_1=i_3$ and $i_2=i_4$, or $i_1=i_4$ and $i_2=i_3$. From
H\"older's inequality we therefore obtain
\[ \Vert g1_V \Vert_A := \Vert (g1_V)^{\wedge} \Vert_1 \geq \frac{\Vert (g1_V)^{\wedge} \Vert_2^3}{\Vert (g1_V)^{\wedge} \Vert_4^2} \geq \big(\frac{|G|}{3}  \big)^{1/2}\Vert (g1_V)^{\wedge} \Vert_2  \geq  \sqrt{m/3}.\]
Since $m > 12M^4$, this is contrary to \eqref{g-l1-bd}, and this
proves the claim.

Applying Lemma \ref{additive-comb-lem} we obtain a subgroup $H'$ such
that
\[ \E_{x \in A}{1_{H'}(x)} \geq 2^{-2^{CM^4}} \; \textrm{ and } \; \E_{x \in {H'}}{1_A(x)} \geq 2^{-CM^4}.\] Since $f^2 \geq 1_A/4$ we get $\E_{x \in
{H'}}f^2(x) \geq 2^{-CM^4}$, but this does
not quite imply Proposition \ref{concentrate-prop}. By
Plancherel's theorem, however, we do have
\[ \langle (f\mu_{H'})^\wedge, \widehat{f}\rangle =\langle f\mu_{H'},f\rangle =  \E_{x \in {H'}}{f^2(x)} \geq 2^{-CM^4},\]
which, since $\Vert f \Vert_A \leq M$, means that
\[ \Vert (f\mu_{H'})^\wedge\Vert_{\infty} \geq 2^{-CM^4}/M \geq 2^{-C'M^4}.\] By the definition of the Fourier transform this yields an $r$ such that
\[ |\E_{x\in {H'}} f(x) (-1)^{r^T x} |\geq 2^{-CM^4}.\]
Taking $H = {H'} \cap \{r\}^{\perp}$, it is clear that
\[ \Vert \psi_{H} f \Vert_{\infty} \geq \sup_{x \in H'} |\E_{y \in x + H} f(y)| \geq 2^{-CM^4}.\]
Finally we note that
\[
\E 1_{H} \geq \E 1_{H'}/2  \geq  2^{-2^{CM^4}} \E 1_A \]
and that
\[ \E1_A \geq \frac{\Vert f_{\mathbb{Z}}\Vert_1}{\|f_\mathbb{Z}\|_\infty}  \geq \frac{\Vert f_{\Z} \Vert_1 }{\|f\|_{\infty} + \epsilon} \geq \frac{\Vert f_{\Z} \Vert_1 }{\|f\|_A + \epsilon} \geq \frac{\|f_\mathbb{Z}\|_1}{M+\epsilon} \geq \frac{\|f_\mathbb{Z}\|_1}{M+1},\] two estimates which together imply the claimed lower bound on $\E 1_H$.\endproof

\section{The main argument}

The basic strategy for proving Theorem \ref{mainthm} is that of an
induction on $M$. Our first lemma provides the main inductive
step. The most noteworthy feature of this lemma is that, in order
to make the induction work, one cannot restrict attention to
boolean functions $f : G \rightarrow \{0,1\}$. It is necessary to
consider almost integer-valued functions as well.

\begin{lemma}[Inductive step]\label{ind-step} Suppose that $f : G
\rightarrow \R$ is $\epsilon$-almost integer-valued with $\Vert f
\Vert_A \leq M$, where $\epsilon \leq 2^{-CM^4}$. Then we may
decompose $f$ as $f_1 + f_2$, where each $f_i$ is
$\epsilon'$-almost integer-valued for some $\epsilon' \leq
2^{CM}\epsilon$, and for $i = 1,2$ one of the following two
alternatives holds:\begin{enumerate} \item $(f_i)_{\Z}$ may be
written as $\sum_{j = 1}^L \pm 1_{H_j}$, where each $H_j$ is a
subgroup of $G$ and
\[ L \leq 2^{2^{CM^4}/\epsilon};\]
\item $\Vert f_i \Vert_A$ is most $\Vert f \Vert_A - \frac{1}{2}$.
\end{enumerate}
\end{lemma}

\proof If $M \leq 1/2$ then $f_{\Z} = 0$, and so option (1)
vacuously holds. Assume, then, that $M \geq 1/2$. We begin by
applying Proposition \ref{concentrate-prop}. This provides a
subgroup $H \leq G$ such that
\[ \E 1_H \geq 2^{-2^{CM^4}} \Vert f_{\Z}\Vert_1 \;\textrm{ and }\; \Vert \psi_H f \Vert_{\infty} \geq
2^{-CM^4}.\]
Set
\[ \eta := 2^{-CM (1+\log M)}\epsilon \] for some large $C$, this choice being dictated by a later application of Proposition \ref{prop3.4}.

By Lemma \ref{lem3.2} we may find a subgroup $H' \leq H$ with
\[ \codim (H : H') \leq M/\eta\]
such that $f$ is $\eta$-spectrally supported on $H'$. By averaging
we have
\begin{equation}\label{eq5.33} \Vert \psi_{H'} f \Vert_{\infty} \geq \Vert \psi_H f \Vert_{\infty} \geq 2^{-CM^4};\end{equation}
we also have
\begin{equation*}
\E 1_{H'} \geq 2^{-M/\eta} \E 1_H \geq
2^{-2^{CM^4}/\epsilon} \|f_{\Z}\|_1.
\end{equation*}

Define $f_1 := \psi_{H'} f$ and $f_2 := f - \psi_{H'} f$. Since
$f$ is $\eta$-spectrally supported on $H'$, it is an immediate
consequence of Proposition \ref{prop3.4} that both $f_1$ and $f_2$
are $\epsilon'$-almost integer-valued, for some $\epsilon' \leq
2^{CM}\epsilon$.

It turns out that for $f_2$ alternative (2) always holds, that is
to say $\Vert f_2 \Vert_A \leq \Vert f \Vert_A - \frac{1}{2}$.
From the Fourier definition of $\psi_{H'}$ one sees that the
supports of $\widehat{f}_1$ and $\widehat{f}_2$ are disjoint, and
hence that
\[ \Vert f \Vert_A = \Vert f_1 \Vert_A + \Vert f_2 \Vert_A.\] Thus we need only show that $\Vert f_1 \Vert_A \geq 1/2$. To see
this, note that from \eqref{eq5.33} we have
\[ \Vert (f_1)_{\Z} \Vert_{\infty} \geq \Vert f_1 \Vert_{\infty} - \epsilon' \geq 2^{-CM^4} - \epsilon' >  0.\]
Since $(f_1)_{\Z}$ is integer-valued, this of course means that
\[ \Vert (f_1)_{\Z} \Vert_{\infty} \geq 1,\]
whence
\[ \|f_1\|_A \geq \Vert f_1 \Vert_{\infty} \geq 1 - \epsilon' \geq 1/2.\]

To conclude the proof, then, we need only show that if $\Vert f_2
\Vert_A < 1/2$ then $(f_1)_{\Z}$ may be written as a $\pm 1$ sum
of not too many cosets of $H'$. The hypothesis on $\|f_2\|_A$ then
ensures that $\|f_2\|_\infty < \epsilon'$ which is certainly at
most $1/10$. Also $f$ is $1/10$-almost integer valued so
\[ \Vert f_{\Z} - \psi_{H'} f \Vert_{\infty} \leq 1/10 + \| f - \psi_{H'}f\|_\infty \leq 1/5.\]
Thus $f_{\Z}$ is within $1/5$ of a function which is constant on
cosets of ${H'}$. Since $f_{\Z}$ is integer-valued, this can only
be the case if $f_{\Z}$ is itself constant on cosets of ${H'}$,
that is to say
\[ f_{\Z} = \sum_{j=1}^{L} c_j 1_{x_j + {H'}}\] for some $x_1,\dots,x_{L}$ which are distinct modulo ${H'}$ and
some non-zero integers $c_j$, $|c_j| \leq \Vert f_{\Z}
\Vert_{\infty} \leq \|f\|_A + \epsilon \leq 2M$.

Recall that the subgroup ${H'}$ is such that $\E 1_{H'} \geq
2^{-2^{CM^4}/\epsilon} \Vert f_{\Z} \Vert_1$, and since we
obviously have
\[ \Vert f_{\Z} \Vert_1 \geq L \E 1_{H'},\] it follows immediately that $L \leq
2^{2^{CM^4}/\epsilon}$. The result follows upon noting that any coset in $G$ is either a subgroup, or else its characteristic function can be expressed as $1_{H_1} - 1_{H_2}$ for two subgroups $H_1,H_2 \leq G$.\endproof

\emph{Proof of Theorem \ref{mainthm}.} We apply Lemma
\ref{ind-step} iteratively, starting with the observation that if
$f : G \rightarrow \{0,1\}$ is a boolean function then $f$ is
$0$-almost integer-valued, and hence an $\epsilon_0$-almost
integer-valued for any $\epsilon_0 > 0$. An appropriate choice of
$\epsilon_0$ will be made later. Split
\[ f = f_1 + f_2\] according to Lemma \ref{ind-step}. Each $f_{i}$ is $\epsilon_1$-almost integer-valued, where
\[ \epsilon_1 \leq 2^{CM}\epsilon_0,\]
and is such that either $(f_i)_{\Z}$ is a sum of at most $2^{2^{CM^4}/\epsilon}$ functions of the form $\pm
1_H$ (in which case we say it is \emph{finished}), or else we have
$\Vert f_i \Vert_A \leq M - \frac{1}{2}$.

Now split any unfinished functions $f_i$ using Lemma
\ref{ind-step} again, and so on (we will discuss the admissibility
of this shortly). This procedure will result in the definition of
parameters $\epsilon_0 \leq \epsilon_1 \leq \dots$ satisfying
\begin{equation}\label{increment}\epsilon_{j+1} \leq 2^{CM}\epsilon_j \end{equation} for all $j$. After at most $2M - 1$ steps all functions will be either finished or will have $\Vert \cdot \Vert_A$-norm at most $1/2$, in which case they are finished for trivial reasons. Thus we have a decomposition
\[ f = \sum_{k = 1}^K f_k,\]
where $K \leq 2^{2M}$, each function $f_k$ is
$\epsilon_{2M}$-almost integer-valued and, for each $k$,
$(f_k)_{\Z}$ may be written as a sum of at most
$2^{2^{CM^4}/\epsilon_0}$ functions $\pm 1_H$.
Now if $\epsilon_0$ is chosen so that $\epsilon_{2M} \leq
2^{-2M}/5$ then
\[ \Vert f - \sum_{k = 1}^K (f_k)_{\Z} \Vert_{\infty} \leq 1/5.\]
However $f$ only takes values in $\{0,1\}$ so it follows that
\[ f = \sum_{k = 1}^K (f_k)_{\Z},\]
which means that $f$ can be written as a sum of at most
$2^{2^{CM^4}/\epsilon_0}$ functions $\pm 1_H$.

The condition $\epsilon_{2M} \leq 2^{-2M}/5$ is not the strongest
condition that we require on $\epsilon_0$. In the repeated
applications of Lemma \ref{ind-step} we must ensure that $\epsilon
\leq 2^{-CM^4}$ is always satisfied, and so we require
\[ \epsilon_{2M} \leq 2^{-CM^4}.\]
In view of \eqref{increment} it is clear, however, that we may
choose $\epsilon_0 \geq 2^{-CM^4}$, for some suitably large $C$,
so that this is indeed always satisfied. This concludes the proof
of Theorem \ref{mainthm}.\endproof

\section{Concluding remarks}

Note that our proof of Theorem \ref{mainthm} actually proves the
following slightly stronger result.
\begin{proposition}
Suppose that $f : G \rightarrow \R$ is a function with $\Vert f
\Vert_A \leq M$, and which is $\epsilon$-almost integer-valued for
some $\epsilon \leq 2^{-CM^4}$. Then the function $f_{\Z}$ can
be written as a combination of at most $2^{2^{CM^4}}$ functions
of the form $\pm 1_H$.
\end{proposition}
In the same way that Theorem \ref{mainthm} can be seen as a
quantitative version of Cohen's result, this proposition can be
seen as a quantitative version of some results of M\'ela
\cite{BHJFMFP,JFM}. An example in M\'ela's work
\cite{BHJFMFP,JFM} shows that $\epsilon$ must be
smaller than $2^{-cM}$ for such a theorem to hold.

\section*{Acknowledgements} The authors are very grateful to Ryan O'Donnell for supplying the
references \cite{MB} and \cite{YM}.

\bibliographystyle{alpha}

\bibliography{master}

\end{document}